\documentclass[12pt]{article}
\usepackage{amssymb,color,todonotes}
\usepackage{thmtools, thm-restate}
\declaretheorem{theorem}

\usepackage[all,knot,arc]{xy}

\newtheorem{prop}[theorem]{Proposition}
\newtheorem{lemma}[theorem]{Lemma}
\newtheorem{cor}[theorem]{Corollary}
\newtheorem{unnumber}{}

\newtheorem{prepf}[unnumber]{Proof}
\newenvironment{pf}{\prepf\rm}{\endprepf}
\newtheorem{preprob}[theorem]{Problem}
\newenvironment{prob}{\preprob\rm}{\endpreprob}

\newcommand{\qed}{\hfill$\Box$}

\newcommand{\Gr}{\mathop{\mathrm{Gr}}}
\newcommand{\pgl}{\mathop{\mathrm{PGL}}}
\newcommand{\End}{\mathop{\mathrm{End}}}
\newcommand{\Aut}{\mathop{\mathrm{Aut}}}
\newcommand{\rank}{\mbox{{rank}}}
\renewcommand{\wr}{\mathbin{\mathrm{wr}}}

\begin{document}
\title{Primitive Groups Synchronize Non-uniform Maps of Extreme Ranks}
\author{Jo\~ao Ara\'ujo\\
  {\small Universidade Aberta, R. Escola Polit\'{e}cnica, 147}\\
  {\small 1269-001 Lisboa, Portugal}\\{\footnotesize \&}\\
  {\small Centro de \'{A}lgebra, Universidade de Lisboa}\\
  {\small 1649-003 Lisboa, Portugal, jaraujo@ptmat.fc.ul.pt}\\\\
  Peter J. Cameron\\
  {\small Mathematical Institute}\\
  {\small North Haugh, St Andrews KY16 9SS, UK}\\
  {\small  pjc@mcs.st-andrews.ac.uk}}
\date{}

\maketitle
\begin{abstract} 
Let $\Omega$ be a set of cardinality $n$, $G$ a permutation group on $\Omega$,
and $f:\Omega\to\Omega$ a map which is not a permutation. We say that $G$
synchronizes $f$ if the semigroup $\langle G,f\rangle$ contains a constant
map.

The first author has conjectured that a primitive group synchronizes
any map whose kernel is non-uniform. Rystsov proved one instance of this conjecture, namely, degree $n$  primitive groups synchronize maps of rank $n-1$ (thus, maps with kernel type $(2,1,\ldots,1)$). We prove some extensions of Rystsov's result, including this: a primitive group synchronizes every map whose kernel type is $(k,1,\ldots,1)$.  Incidentally this result provides a new characterization of imprimitive groups.  We also  prove that  the conjecture above holds for maps of extreme ranks, that is, ranks 3, 4 and $n-2$. 
 
 These proofs use a graph-theoretic
technique due to the second author: a transformation semigroup fails to
contain a constant map if and only if it is contained in the endomorphism
semigroup of a non-null (simple undirected) graph. 

The paper finishes with a number of open problems, whose solutions will certainly require very delicate graph theoretical considerations.
\end{abstract}

\section{Introduction}
In automata theory, the well-known  \v Cern\'y conjecture states that a synchronizing automaton with $n$ states has a synchronizing word of length $(n-1)^2$. (For many references on the growing bibliography on this problem please see the two websites \cite{JEP,Tr}.)
Solving this conjecture is equivalent to prove that given a set $S=\{f_{1},\ldots,f_{m}\}$ of transformations on a finite set $\Omega:=\{1,\ldots,n\}$, if $S$ generates a constant, then $S$ generates a constant in a length $(n-1)^2$ word on its generators.  
This conjecture has been established when $\langle S\rangle$  is a semigroup in which all its subgroups are trivial  \cite{Tr07}. So it remains to prove the
conjecture for semigroups containing non trivial subgroups; the case
in which the semigroup contains a permutation group is  a particular instance of this
general problem.  In addition, the known examples  witnessing the optimality of the \v Cern\'y bound contain a permutation among the given set of generators $S$, so they
make it especially interesting to study the cases in which a subset of $S$ generates a permutation group.

Let $G$ be a permutation group on a set $\Omega$ with $|\Omega|=n$.
We say that $G$ \emph{synchronizes} a map $f$ on $\Omega$ if the semigroup
$\langle G,f\rangle$ contains a constant map. $G$ is said to be \emph{synchronizing} if $G$ synchronizes every non-invertible transformation on $\Omega$. The \emph{diameter} of a group is the largest diameter of its Cayley graphs. Taking into account the motivation of the considerations above, the ultimate goal is to find a classification of the synchronizing groups and then study those with the largest diameter, since they should assist the generation of a constant with the lowest diligence. But even when we forget about the automata motivation of these problems, the
  classification of synchronizing groups (a class strictly between primitivity and $2$-homogeneity) and the study of their diameters are very interesting questions in themselves, as well as extremely  demanding (please see \cite{abc,ArnoldSteinberg,bamberg,neumann,helfgott}). 

Let $f$ be a map on $\Omega$. Recall that the \emph{rank} of  $f$  is $|\Omega f|$,
and the \emph{kernel} of $f$  is
the partition of $\Omega$ into the inverse images of points in the image
of $f$; equivalently, the kernel of $f$ is the partition of $\Omega$ induced by the equivalence relation $\{(x,y)\in \Omega\times \Omega \mid xf=yf\}$. The \emph{kernel type} of $f$ is the partition of $n$ given by the
sizes of the parts of the kernel. A partition of $\Omega$ is \emph{uniform}
if all its parts have the same size. We will call a map \emph{uniform}
if its kernel is uniform.

 We note that, if a transformation
semigroup $S$ contains a transitive group $G$ but not a constant function, then
the image $I$ of a map $f$ of minimal rank in $S$ is a
\emph{$G$-section} for the kernel of $f$, in the sense that $Ig$ is a section for $\ker (f)$, for
all $g\in G$; in addition, the map $f$ has uniform kernel  (see Neumann~\cite{neumann}). 

In \cite{abc}  the conjecture that a primitive group of permutations of $\Omega$ synchronizes every non-uniform transformation on $\Omega$ was proposed.
In 1995 Rystsov~\cite{rystsov} proved the following particular instance of this conjecture. 

\begin{theorem}\label{rys}
A transitive permutation group $G$ of degree $n$ is primitive if and only if
it synchronizes every map of rank $n-1$.
\end{theorem}

The goal of this paper is to use a graph-theoretic approach due to the second author  to prove the conjecture  for maps of extreme rank, that is, their rank either is close to $1$ or  close to $n$.

It is worth pointing out that Rystosov's Theorem \ref{rys} in fact characterizes primitivity in terms of maps of kernel type $(2,1,\ldots,1)$; our first result provides a similar characterization of  imprimitivity in terms of maps of kernel type $(k,1,\ldots,1)$.

\begin{restatable}{theorem}{imprimitive}
\label{zerost}
Let $G$ be a transitive permutation group on a finite set $\Omega$ and let 
 $k$ be given with $k>1$. Then $G$ is imprimitive with a block of
imprimitivity of size at least $k$ if and only if $G$ fails to synchronize
some map $f$ with kernel type $(k,1,1,\ldots,1)$.
\end{restatable}

This result proves, in particular, that primitive groups synchronize every map with kernel type $(k,1,1,\ldots,1)$.

Our next result concerns synchronization of large rank maps.

\begin{restatable}{theorem}{primitive}
\label{first}
Let $G$ be a primitive permutation group on a finite $\Omega$, with $|\Omega|>2$. Then $G$ synchronizes:
\begin{enumerate}
\item\label{b}  every map of rank $n-2$;
\item\label{c}   every idempotent map with
kernel type $(3,2,1,1,\ldots,1)$;
\item\label{d}   every map $f$ with
kernel type $(3,2,1,1,\ldots,1)$, provided  there exists $g\in G$ such that  $\rank(fgf)=\rank(f)$.
\end{enumerate}
\end{restatable}

In the second part of the paper, we turn from maps of large rank to those
of small rank, and prove the following theorem. (The first part of this
theorem is due to Neumann~\cite{neumann}, but we will provide an alternative proof.)

\begin{restatable}{theorem}{small}
\label{second}
Let $G$ be a primitive group of degree $n>2$. 
\begin{enumerate}
\item $G$ synchronizes every map of rank $2$.
\item $G$ synchronizes every non-uniform map of rank $3$ or $4$.
\end{enumerate}
\end{restatable}


The condition that the map is non-uniform in Theorem~\ref{second}(b) is necessary:
the group $S_3\wr S_2$ of degree~$9$ (the automorphism group of the $3\times3$
grid) is primitive and fails to synchronize a map of rank $3$ (for example,
the projection of the grid onto a diagonal whose kernel classes are the rows).

In Section \ref{trans} we introduce the graph, and some basic results about it, that is going to be our main tool throughout the paper. Section \ref{impri} is dedicated to the proof of Theorem \ref{zerost}, Section \ref{pri} to the proof of Theorem \ref{first}, and Section \ref{sma} to the proof of Theorem \ref{second}. The paper ends with a number of open problems whose solution will certainly require delicate considerations on graph theory, permutations groups and transformation semigroups. 

\section{Transformation semigroups and graphs}\label{trans}

The critical idea used in this paper is a graph associated to a transformation semigroup,  due to the second author. Let $S$ be a transformation semigroup on $\Omega$. Form a graph on the
vertex set $\Omega$ by joining two vertices $v$ and $w$ if and only if there
is no element $f$ of $S$ which maps $v$ and $w$ to the same point. We denote
this graph by $\Gr(S)$. Now the following result is almost immediate (\emph{cf.} \cite{CK}).
\begin{theorem} Let $S$ be a transformation semigroup on $\Omega$ and let $\Gr(S)$ be as above. 

\begin{enumerate}
\item[(a)] $S$ contains a map of rank $1$ if and only if $\Gr(S)$ is null.
\item[(b)] $S\le\End(\Gr(S))$, and $\Gr(\End(\Gr(S)))=\Gr(S)$.
\item[(c)] The clique number and chromatic number of $\Gr(S)$ are both equal to
the minimum rank of an element of $S$.
\end{enumerate}
\end{theorem}
\begin{pf}
Regarding (a), the forward direction is obvious. Conversely, let $f\in S$ be a map of minimal rank, and suppose that $\rank(f)>1$. For every $x,y\in \Omega f$ we have $xs\neq ys$, for all $s\in S$ (otherwise $\rank(fs)<\rank(f)$ contrary to our assumption).   Therefore $\{x,y\}$ is an edge of $\Gr(S)$. It is proved that if $S$ has no constant, then $\Gr(S)$ is not null. 

Regarding (b), let $f\in S$ and let $\{x,y\}$ be any edge in $\Gr(S)$; we claim that $\{xf,yf\}$ is an edge in $\Gr(S)$ and hence $f\in \End(\Gr(S))$. In fact, if $\{xf,yf\}$ is not an edge in $\Gr(S)$, then there exists $f'\in S$ such that $xff'=yff'$, that is,  $\{x,y\}$ is not an edge of $\Gr(S)$, contradicting our assumptions. 

Now we prove that $\Gr(\End(\Gr(S)))\subseteq\Gr(S)$. Let $\{x,y\}$ be an edge in $\Gr(\End(\Gr(S)))$. This means that $xf\neq yf$, for all $f\in \End(\Gr(S))$. But we already proved that  $S\le\End(\Gr(S))$; thus $xf\neq yf$, for all $f\in S$ and hence $\{x,y\}$ is an edge in $\Gr(S)$. 

Conversely, to prove that $\Gr(\End(\Gr(S)))\supseteq\Gr(S)$, let $\{x,y\}$ be an edge in $\Gr(S)$. Then for every $f\in \End(\Gr(S))$ the set $\{xf,yf\}$ is an edge of $\Gr(S)$, that is, $xf\neq yf$, for all $f\in \End(\Gr(S))$. Thus $\{x,y\}$ is an edge of $\Gr(\End(\Gr(S)))$ and (b) follows. 

Now we prove (c). It is clear that the image of any map $f$ of minimum rank forms a clique of $\Gr(S)$; for if not there would be $x,y\in \Omega f$ and $f'\in S$ such that $xf'=yf'$; thus $\rank(ff')$ would be strictly smaller than $\rank (f)$ and hence $f$ would not be of minimum rank. Let $\Gamma$ be the complete graph contained in $\Gr(S)$ and whose vertex set is $\Omega f$ (for a map $f$ of minimum rank). It is clear that $f:\Gr(S)\mapsto \Gamma$ is a morphism; conversely, $\iota: \Gamma \mapsto \Gr(S)$ such that $\{x,y\}\iota = \{x,y\}$ is a morphism. Thus the complete graph $\Gamma$ is a core of $\Gr(S)$; it is well known that if a graph has complete core, then the chromatic number of the graph equals its clique number. The result follows.    

Note that (a) is a special case of (c), when the minimum rank is~$1$.
\qed
\end{pf}

In particular, if $S=\langle G,f\rangle$ for some group $G$, then
$G\le\Aut(\Gr(S))$. So, for example, if $G$ is primitive and does not
synchronize $f$, then $\Gr(S)$ is non-null and has a primitive automorphism group, and so is
connected.

In this situation, assume that $f$ is an element of minimal rank in $S$;
then the kernel of $f$ is a partition $\rho$ of $\Omega$, and its image $A$
is a \emph{$G$-section} for $\rho$ (that is, $Ag$ is a section for $\rho$,
for all $g\in G$). Neumann~\cite{neumann}, analysing this situation, defined
a graph $\Delta$ on $\Omega$ whose edges are the images under $G$ of the
pairs of vertices in the same $\rho$-class. Clearly $\Delta$ is a subgraph
of the complement of $\Gr(S)$, since edges in $\Delta$ can be collapsed by
elements of $S$. Sometimes, but not always, $\Delta$ is the complement of
$\Gr(S)$.

For the sake of completeness we include here a general lemma on primitive groups.  
\begin{lemma}\label{transp}
Let $G$ be  a primitive group:
\begin{enumerate}
\item if $G$ contains a transposition $(v,w)$, then $G$ is the symmetric group;
\item if $G$ has degree greater than $5$ and contains a double transposition $(v,w)(x,y)$, then $G$ is $2$-transitive.
\end{enumerate}
 \end{lemma}
\begin{pf}
Regarding (a), suppose $G$ is a primitive group of permutations of $\Omega:=\{1,\ldots ,n\}$. Define a relation on $\Omega$ as follows: for all $x,y\in \Omega$,
\[
x\sim y \Leftrightarrow x=y \mbox{ or } (x,y) \in G. 
\]

It is clear that $\sim$ is reflexive and symmetric. In addition, if $x\sim y\sim z$, then $(x,y),(y,z)\in G$ and hence $(x,z) =(y,z)(x,y)(y,z)\in G$; thus $\sim$ is an equivalence relation on $\Omega$. The transpositions generate the symmetric groups on the equivalence classes. We claim that there is only one equivalence class; for suppose not and let $(x,y)\in G$ and $A\subset \Omega$ be an equivalence class. Since $G$ is primitive there exists $g\in G$ such that $xg \in A$ and $yg \not\in A$. Thus $g^{-1}(x,y)g=(xg,yg)\in G$; thus $xg\sim yg$, a contradiction. It is proved that there is only one equivalence class and it was already shown above that the transpositions generate  the symmetric group inside each equivalence class.

Regarding (b), we refer the reader to Example 3.3.1 on p.82 of \cite{DM}. 
\qed
\end{pf}
The next lemma has some interest in itself, but it is very important for the techniques it introduces and that will be used later.

\begin{lemma}\label{neigh}
Let $X$ be a nontrivial graph and let $G\le \Aut(X)$ be primitive. Then no two vertices of $X$ can have the same neighbourhood. 
\end{lemma}
\begin{pf}
For $a\in X$ denote its neighbourhood by $N(a)$. Suppose that $a,b\in X$, with $a\neq b$, and $N(a)=N(b)$. We are going to use two different techniques to prove that this leads to a contradiction. The first uses the fact that the graph has at least one edge; the second uses the fact that the graph is not complete. 

First technique. Define the following relation on the vertices of the graph: for all $x,y\in X$,
\[
x\equiv y \Leftrightarrow N(x)=N(y). 
\] 
This is an equivalence relation and we claim that $\equiv$ is neither the universal relation nor the identity. The latter follows from the fact that by assumption $a$ and $b$ are different and $N(a)=N(b)$. Regarding the former, there exist adjacent vertices $c$ and $d$ (because $X$ is non-null); 
now $c \in N(d)$
but $c \notin N(c)$, so $c \not\equiv d$.
As $G$ is a group of automorphisms of $X$ it follows that $G$ preserves $\equiv$, a non-trivial equivalence relation, and hence $G$ is imprimitive, a contradiction.

Second technique. Assume as above that we have $a,b\in X$ such that $N(a)=N(b)$. Then the transposition $(a,b)$ is an automorphism of the graph. By the previous lemma, a primitive group containing a transposition is the symmetric group and hence $X$ is the complete graph, a contradiction.  
\qed
\end{pf}

The two techniques in the previous proof are important because we will use
variants on them later.

We conclude this section with a general result about \emph{primitive graphs}
(those admitting a vertex-primitive automorphism group), which we will use
later in the paper.

\begin{lemma}
Let $\Gamma$ be a non-null graph with primitive automorphism group $G$, and
having chromatic number $r$. Then $\Gamma$ does not contain a subgraph 
isomorphic to the complete graph on $r+1$ vertices with an edge removed.
\label{primgr}
\end{lemma}
\begin{pf}
Let $c$ be a colouring of $\Gamma$ with $r$ colours. Suppose, by contradiction, that the set
$\{1,2,\ldots,r,r+1\}$ of vertices contains all possible edges except for
$\{r,r+1\}$. Then $\{1,\ldots,r\}$ is a clique, and so contains one vertex
of each colour; similarly for $\{1,\ldots,r-1,r+1\}$. Since the colors of $r$ and $r+1$ are different from the colors of $1,\ldots,r-1$, we conclude that
vertices $r$ and $r+1$ have the same colour. The same conclusion holds for
the image of these vertices under any element of $G$.

Now let $\Delta$ be the graph whose edge set is the $G$-orbit containing
$\{r,r+1\}$. Then $\Delta$ is $G$-invariant and non-empty, but is disconnected,
since all its edges lie within colour classes of the colouring $c$. This
contradicts the primitivity of $G$.
\qed
\end{pf}

We note that the hypotheses are both necessary: the complete $r$-partite graph
with parts of constant size is vertex-transitive and contains $K_{r+1}$ minus
an edge; and every graph occurs as an induced subgraph of some primitive graph  as proved in the next result.

\begin{prop}
Every graph is isomorphic to an induced subgraph of a graph with primitive
automorphism group.
\end{prop}
\begin{pf}
First represent the graph as an \emph{intersection graph}, that is, the
vertices are subsets of a set $E$, and two vertices are adjacent if the sets
are not disjoint. This was first observed by Szpilrajn-Marczewski \cite{S-M};
it is most easily done by taking $E$ to be the edge set of the graph, and
identifying each vertex with the set of edges incident with it.

Now, by adding extra points each in at most one of the sets, we may
assume that all the sets have the same cardinality $k$.

Now the graph is an induced subgraph of the graph whose vertices are the
$k$-element subsets of an $n$-set (where we may assume that $n>2k$), two
vertices adjacent if they are not disjoint. The automorphism group of this
graph is the symmetric group $S_n$, in its primitive action on $k$-sets.\qed
\end{pf}

\section{A characterization of imprimitivity}\label{impri}

Rystsov's Theorem  says that every primitive group synchronizes a map of kernel type $(2,1,\ldots,1)$. The next theorem generalizes this result by proving that every primitive group synchronizes a map of kernel type $(k,1,\ldots,1)$, for every $k$ such that $|\Omega|\geq k\geq 2$.

In fact, using the graph-theoretic techniques of the preceding section we prove the following characterization of imprimitivity.


\imprimitive*

\begin{pf} To prove the theorem in the forward direction,
 suppose that $G$ is imprimitive, with blocks of size at
least $k$. Let $X$ be the complete multipartite graph whose partite classes
are the blocks. Then $G\le\Aut(X)$. Let $A$ be a subset of a block,
with $|A|=k$, and choose $a\in A$. Define $f$ so that $bf=a$ for all $b\in A$
and $xf=x$ for $x\notin A$. Then $f$ is an endomorphism of $X$ (so that
$f$ cannot be synchronized by $G$) with kernel type $(k,1,1,\ldots,1)$.

Conversely,  let $G$ and $f$ be as given, and
let $A$ be the kernel class of size $k$ of $G$. 

By the general observations of the previous section  we know that there is a non-null graph
$X=\Gr(\langle G,f\rangle)$ with $\langle G,f\rangle\le\End(X)$.

Observe that $A$ is an independent set in $X$, since $A$ is collapsed
to a point by the endomorphism $f$. Thus $X$ is not the complete graph and hence, since $X$ is also not null, we conclude that $X$ is not trivial. 

We claim that any two points in $A$ have the same neighbourhood. To see that let
$N(x)$ denote the neighbourhood of $x$. Since $G\leq \Aut (X)$ is
transitive, all vertices have the same number of neighbours. Now let
$Af=\{z\}$. Then, for any $a\in A$, as $f\in \End(X)$ it follows that $f$ maps $N(a)$ to $N(z)$; since
$N(a)\cap A=\emptyset$, $f$ is injective on $N(a)$, and so maps it bijectively
to $N(z)$ (recall from  above that  the transitivity of $G$ implies that   $|N(a)|=|N(z)|$). Similarly for another point $b\in A$. But $z\notin N(z)$, so
$f^{-1}|_{{N(z)}}$ is a  well defined bijective map from $N(z)$ to $N(a)$, and also
to $N(b)$; so these two sets are equal. By Lemma \ref{neigh} this immediately implies that $G$ is imprimitive.

Now recall the $G$-congruence $\equiv$ on $\Omega$ defined in the proof of Lemma \ref{neigh}:  $x\equiv y$ if and
only if $N(x)=N(y)$. Since all elements in $A$ have the same neighbourhood, we conclude that $A$ is contained in a single $\equiv$-class.   So $G$ is
imprimitive, with a system of blocks (the $\equiv$-classes) of size at
least $k$.\qed
\end{pf}

\section{Primitive groups and large rank maps}\label{pri}
The aim of this section is to prove the following theorem. 

\primitive*

The proof of this theorem will be carried out in a sequence of subsections.
\subsection{Proof of Theorem \ref{first}(\ref{b})}

The kernel type of a map of rank $n-2$ is either $(3,1,1,\ldots,1)$ or
$(2,2,1,\ldots,1)$. By Theorem \ref{zerost}, a primitive group synchronizes a map
with the first kernel type. 

We begin with some general remarks about the case where $f$ is a map with
exactly two non-singleton kernel classes $A$ and $B$, and $G$ is a primitive group
which fails to synchronize $f$. 

We start with a general lemma.  

\begin{lemma}\label{quadrangle}
Given a primitive group $G$ that does not synchronize a map of kernel type $(p,q,1,\ldots,1)$, with $p,q\ge 2$. Let $A$ and $B$ be the non-singleton kernel classes, and $K:=A\cup B$. Let
$S=\langle G,f \rangle$ and $X=G(S)$. Then   there must be at least a path or a cycle of length $4$ contained in the non-singleton kernel classes $A$ and $B$. Moreover, let $S=\langle G,f\rangle$ and $X=\Gr(S)$. Then in   the graph $X$ we have that 
\begin{enumerate}
\item there are no isolated points in $K$;
\item there cannot be two points of $A$ which each have a single
neighbour in $B$, these neighbours being equal.
\end{enumerate}
\end{lemma}
\begin{pf}
The rank of $f$ must be larger than $2$ since, by \cite{neumann} and/or Theorem \ref{second},  $G$ synchronizes every rank $2$ map; therefore, $A\cup B$ is a proper subset of $\Omega$. We let $S=\langle G,f\rangle$ and $X=\Gr(S)$.
Suppose that $Af=x$ and $Bf=y$.

\begin{figure}[h]
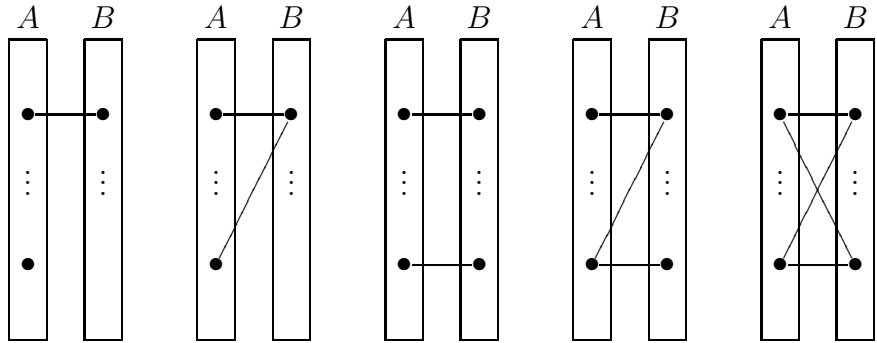

\[
\xy
(-25,0)*{}="0a";
(-25,40)*{}="0b";
(-20,40)*{}="0c";
(-20,0)*{}="0d";
(-22.5,22)*{\vdots};
(-22.5,10)*{\bullet}="0e";
(-22.5,30)*{\bullet}="0f";
(-22.5,43)*{A};
"0a";"0b" **\dir{-};
"0b";"0c" **\dir{-};
"0c";"0d" **\dir{-};
"0a";"0d" **\dir{-};
(-15,0)*{}="01a";
(-15,40)*{}="01b";
(-10,40)*{}="01c";
(-10,0)*{}="01d";
(-12.5,22)*{\vdots};
(-12.5,10)*{}="01e";
(-12.5,30)*{\bullet}="01f";
(-12.5,43)*{B};
"01a";"01b" **\dir{-};
"01b";"01c" **\dir{-};
"01c";"01d" **\dir{-};
"01a";"01d" **\dir{-};
"0f";"01f" **\dir{-};
(0,0)*{}="a";
(0,40)*{}="b";
(5,40)*{}="c";
(5,0)*{}="d";
(2.5,22)*{\vdots};
(2.5,10)*{\bullet}="e";
(2.5,30)*{\bullet}="f";
(2.5,43)*{A};
"a";"b" **\dir{-};
"b";"c" **\dir{-};
"c";"d" **\dir{-};
"a";"d" **\dir{-};
(10,0)*{}="1a";
(10,40)*{}="1b";
(15,40)*{}="1c";
(15,0)*{}="1d";
(12.5,22)*{\vdots};
(12.5,10)*{}="1e";
(12.5,30)*{\bullet}="1f";
(12.5,43)*{B};
"1a";"1b" **\dir{-};
"1b";"1c" **\dir{-};
"1c";"1d" **\dir{-};
"1a";"1d" **\dir{-};
"e";"1f" **\dir{-};
"f";"1f" **\dir{-};
(25,0)*{}="2a";
(25,40)*{}="2b";
(30,40)*{}="2c";
(30,0)*{}="2d";
(27.5,22)*{\vdots};
(27.5,10)*{\bullet}="2e";
(27.5,30)*{\bullet}="2f";
(27.5,43)*{A};
"2a";"2b" **\dir{-};
"2b";"2c" **\dir{-};
"2c";"2d" **\dir{-};
"2a";"2d" **\dir{-};
(35,0)*{}="3a";
(35,40)*{}="3b";
(40,40)*{}="3c";
(40,0)*{}="3d";
(37.5,22)*{\vdots};
(37.5,10)*{\bullet}="3e";
(37.5,30)*{\bullet}="3f";
(37.5,43)*{B};
"3a";"3b" **\dir{-};
"3b";"3c" **\dir{-};
"3c";"3d" **\dir{-};
"3a";"3d" **\dir{-};
"2e";"3e" **\dir{-};
"2f";"3f" **\dir{-};
(50,0)*{}="4a";
(50,40)*{}="4b";
(55,40)*{}="4c";
(55,0)*{}="4d";
(52.5,22)*{\vdots};
(52.5,10)*{\bullet}="4e";
(52.5,30)*{\bullet}="4f";
(52.5,43)*{A};
"4a";"4b" **\dir{-};
"4b";"4c" **\dir{-};
"4c";"4d" **\dir{-};
"4a";"4d" **\dir{-};
(60,0)*{}="5a";
(60,40)*{}="5b";
(65,40)*{}="5c";
(65,0)*{}="5d";
(62.5,22)*{\vdots};
(62.5,10)*{\bullet}="5e";
(62.5,30)*{\bullet}="5f";
(62.5,43)*{B};
"5a";"5b" **\dir{-};
"5b";"5c" **\dir{-};
"5c";"5d" **\dir{-};
"5a";"5d" **\dir{-};
"4e";"5f" **\dir{-};
"4f";"5f" **\dir{-};
"4e";"5e" **\dir{-};
(75,0)*{}="6a";
(75,40)*{}="6b";
(80,40)*{}="6c";
(80,0)*{}="6d";
(77.5,22)*{\vdots};
(77.5,10)*{\bullet}="6e";
(77.5,30)*{\bullet}="6f";
(77.5,43)*{A};
"6a";"6b" **\dir{-};
"6b";"6c" **\dir{-};
"6c";"6d" **\dir{-};
"6a";"6d" **\dir{-};
(85,0)*{}="7a";
(85,40)*{}="7b";
(90,40)*{}="7c";
(90,0)*{}="7d";
(87.5,22)*{\vdots};
(87.5,10)*{\bullet}="7e";
(87.5,30)*{\bullet}="7f";
(87.5,43)*{B};
"7a";"7b" **\dir{-};
"7b";"7c" **\dir{-};
"7c";"7d" **\dir{-};
"7a";"7d" **\dir{-};
"6e";"7f" **\dir{-};
"6f";"7f" **\dir{-};
"6e";"7e" **\dir{-};
"6f";"7e" **\dir{-};
\endxy
\]
\caption{The five possible configurations of the edges in $K$}\label{f38}
\end{figure}

The graph $X$ is non-trivial and has the primitive group $G$ contained in $\Aut(X)$; so $X$ is
connected, and hence there exists at least one edge from $A\cup B$ to its complement (recall that $A\cup B$ is a proper subset of $\Omega$).

Certainly neither  $A$ nor $B$ contain edges of $X$, since
$f$ collapses each one of them. So $X$ is not complete.

We claim that the set $K:=A\cup B$ must contain an edge of $X$. For suppose not. Take
$a_1,a_2\in A$. Then $N(a_1)\subseteq X\setminus (A\cup B)$ and $f$ is injective on $X\setminus (A\cup B)$. Thus $f$ maps $N(a_1)$ injectively into $N(x)$. The transitivity of $G$ implies that all neighbourhoods have the same size and hence $N(a_1)f=N(x)$; so does $N(a_2)$. Thus $N(a_1)=N(x)f^{-1}=N(a_2)$, which is impossible in a
primitive group (Lemma \ref{neigh}).
 It follows that $K$ contains an edge and hence $\{x,y\}$ is an edge (since $f$ is an
endomorphism of $X$). 
Therefore the edges in $K$ must contain one of the five configurations in figure \ref{f38}. We are going to show that only the last two configurations can happen.

First, we claim that any point of $K$ lies on an edge within
this set. For suppose that $a\in A$ does not; then all its neighbours are
outside $K$, and so are mapped bijectively by $f$; but they are mapped
onto $N(x)\setminus\{y\}$, which is smaller than $N(a)$.
It is proved that the first configuration in figure \ref{f38} cannot occur. 

Next we claim that there cannot be two points of $A$ which each have a single
neighbour in $B$, these neighbours being equal. For suppose that $a_1$ and
$a_2$ are two such points. Then the sets
$N(a_1)\setminus B$ and $N(a_2)\setminus B$ are mapped injectively,
and hence bijectively, to $N(x)\setminus\{y\}$; so we must have
$N(a_1)\setminus B=N(a_2)\setminus B$. But $a_1$ and $a_2$ have the
same neighbour in $B$; so their neighbourhoods are equal, which again is
impossible in a primitive group (Lemma \ref{neigh}). It is proved that the second configuration in figure \ref{f38} cannot occur. 

Also, the induced subgraph on $K$ cannot have two connected components each
consisting of a single edge. For suppose that $\{a_1,b_1\}$ and
$\{a_2,b_2\}$ were such components. As above, we would have
$N(a_1)\setminus B=N(a_2)\setminus B$ and, similarly,
$N(b_1)\setminus A=N(b_2)\setminus A$. But then the permutation
$(a_1,a_2)(b_1,b_2)$ is an automorphism of the graph. However, a primitive
group of degree greater than $5$ containing such an element must be
$2$-transitive (by Lemma \ref{transp}) whereas a non-trivial
graph cannot have a $2$-transitive automorphism group.\qed
\end{pf}

Now we can prove Theorem \ref{first}(\ref{b}). 
The kernel type of a map of rank $n-2$ is either $(3,1,1,\ldots,1)$ or
$(2,2,1,\ldots,1)$. By Theorem \ref{zerost}, a primitive group synchronizes a map
with the first kernel type. 

Let $f$ be a map of kernel type $(2,2,1,\ldots,1)$ and let $S=\langle G,f\rangle$ and $X=\Gr(S)$.
Suppose that $Af=x$ and $Bf=y$.  By the previous result we already know that $K:=A\cup B$ must contain a path or a cycle of length $4$. In order to  finish  the proof we must  consider those two  cases for the induced subgraph on $K$: a path or 
cycle of length $4$, when $|A|=2=|B|$. 

Since $G\leq\Aut(X)$ is transitive, $X$ is regular of valency, say, $k$. Observe that in each case, the $k-1$ or $k-2$ vertices of
$N(a_i)\setminus K$ are mapped injectively to the vertices of
$N(x)\setminus\{y\}$, and the $k-1$ or $k-2$ vertices of
$N(b_i)\setminus K$ to the vertices of $N(y)\setminus\{x\}$.

Now the graph $X f$ is a subgraph of the induced
subgraph of $X$ on $\Omega f$, call it $X'$, say.  Then in $X'$ we have removed two vertices, which are incident
with either $2k-1$ or $2k$ edges, according as they are adjacent or not. So
$X'$ has at most $e-2k+1$ edges, where $e$ is the number of edges of
$X$. We show that this is incompatible with each of our cases except in
one situation.

To get a lower bound for the number of edges of $X f$, we simply have to
calculate the size of the image of the edge set of $X$ under $f$. (Each
edge of $X$ maps to an edge of $X f$.) We do this by counting edges
collapsed by $f$.

Consider the case where $K$ induces a path of length $3$, say
$(a_1,b_1,a_2,b_2)$. Now $f$ collapses the three edges within
$K$ to the single edge $\{x,y\}$. It maps the $k-1$ vertices of
$N(a_1)\setminus K$ bijectively to $N(x)\setminus\{y\}$, and the
$k-2$ vertices of $N(a_2)\setminus K$ injectively inside this set;
so $k-2$ pairs of edges of this form collapse. Similarly $k-2$ pairs of
edges through $b_i$ collapse. So the number of edges of
$X f$ is at least $e-(2+2(k-2))=e-2k+2$, a contradiction.

Now consider the case where $K$ induces a $4$-cycle. The neighbours of
both $a_1$ and $a_2$ outside $K$ are mapped injectively to
$N(x)\setminus\{y\}$, so either $k-2$ or $k-3$ edges are collapsed,
depending on whether the images of these two sets are equal or not.
Similarly for $b_1$ and $b_2$. Moreover, four edges within $K$ are collapsed
to one. So the number of edges of $X f$ is at least $e-3-2(k-2)=e-2k+1$.
This is just possible, but we see that the sets $N(a_1)\setminus K$
and $N(a_2)\setminus K$ must be equal, and similarly
$N(b_1)\setminus K=N(b_2)\setminus K$.

Now, if we assume that $A=\{a_{1},a_{2}\}$ and $B=\{b_{1},b_{2}\}$, the equality of the neighbourhoods just proved shows that
the permutation $(a_1,a_2)(b_1,b_2)$ is an automorphism of $X$, a
contradiction as before. 
This finishes the proof of Theorem \ref{first}(\ref{b}).\qed

\medskip

\subsection{Proof of Theorem \ref{first}(\ref{c})}

As in the previous subsection, we begin with a more general result. Let $f$ be a map with kernel type
$(p,q,1,\ldots,1)$ (with $p,q>1)$ non synchronized by a primitive group $G$; let
$S=\langle G,f\rangle$ and $X=\Gr(S)$. Suppose, further, that $f$ is an
idempotent. We claim that the induced subgraph on $K=A\cup B$ (where $A$ and
$B$ are the non-trivial kernel classes) cannot have all possible edges between
$A$ and $B$.

As in the preceding section, we use the fact that $Xf$ is a subgraph of the
restriction $X'$ of $X$ to the image of $f$. Now the fact that $f$ is an
idempotent means that the points not in the image of $f$ must consist of
$p-1$ points of $A$ and $q-1$ points of $B$; by assumption, the induced
subgraph on this set is complete bipartite.

Let $e$ be the number of edges, and $k$ the valency, of $X$.

In the subgraph $X'$, we lose $p-1$ vertices of $A$ and $q-1$ of $B$,
and so $k(p-q+2)$ edges; but of these, $(p-1)(q-1)$ are counted twice. So the
number of edges in this graph is
\begin{eqnarray}\label{cardinal}
e-k(p+q-2)+(p-1)(q-1)
\end{eqnarray}

Now we consider $Xf$, and count how many edges collapse to the same
place. The $pq$ edges within $K$ collapse to a single edge, so we lose $pq-1$
edges. The $p$ vertices of $A$ each lie on $k-q$ edges outside $K$, so at
most $(p-1)(k-q)$ edges are lost; and similarly at most $(q-1)(k-p)$ edges
through vertices in $B$. Thus the number of edges in $X f$ is at least
\[e-(pq-1)-(p-1)(k-1)-(q-1)(k-p)=e-(p+q-2)k-(p-1)(q-1),\]
and this number equals the value found in (\ref{cardinal}). So all the neighbours of $A$ outside $K$ are
mapped to the same set of $k-q$ vertices. This means that any two vertices
in $A$ have the same neighbours. Therefore any transposition of two elements in $A$ is an automorphism of $X$ and hence $G$ is the symmetric group (by Lemma \ref{transp});  but this is impossible, since $G$ does not synchronize $f$ while the symmetric group synchronizes every map. 

\medskip

\begin{figure}[h]
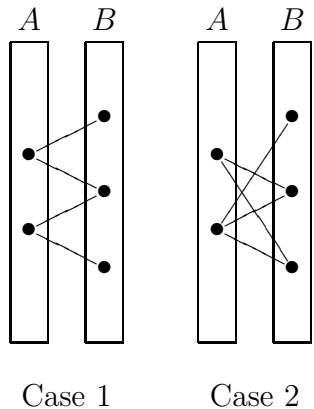

\[
\xy
(7.5,-7)*{\mbox{Case 1}};
(0,0)*{}="a";
(0,40)*{}="b";
(5,40)*{}="c";
(5,0)*{}="d";
(2.5,15)*{\bullet}="e";
(2.5,25)*{\bullet}="f";
(2.5,43)*{A};
"a";"b" **\dir{-};
"b";"c" **\dir{-};
"c";"d" **\dir{-};
"a";"d" **\dir{-};
(10,0)*{}="1a";
(10,40)*{}="1b";
(15,40)*{}="1c";
(15,0)*{}="1d";
(12.5,10)*{\bullet}="1e";
(12.5,20)*{\bullet}="1f";
(12.5,30)*{\bullet}="1g";
(12.5,43)*{B};
"1a";"1b" **\dir{-};
"1b";"1c" **\dir{-};
"1c";"1d" **\dir{-};
"1a";"1d" **\dir{-};
"e";"1f" **\dir{-};
"e";"1e" **\dir{-};
"f";"1f" **\dir{-};
"f";"1g" **\dir{-};
(32.5,-7)*{\mbox{Case 2}};
(25,0)*{}="2a";
(25,40)*{}="2b";
(30,40)*{}="2c";
(30,0)*{}="2d";
(27.5,15)*{\bullet}="2e";
(27.5,25)*{\bullet}="2f";
(27.5,43)*{A};
"2a";"2b" **\dir{-};
"2b";"2c" **\dir{-};
"2c";"2d" **\dir{-};
"2a";"2d" **\dir{-};
(35,0)*{}="3a";
(35,40)*{}="3b";
(40,40)*{}="3c";
(40,0)*{}="3d";
(37.5,10)*{\bullet}="3e";
(37.5,20)*{\bullet}="3f";
(37.5,30)*{\bullet}="3g";
(37.5,43)*{B};
"3a";"3b" **\dir{-};
"3b";"3c" **\dir{-};
"3c";"3d" **\dir{-};
"3a";"3d" **\dir{-};
"2e";"3e" **\dir{-};
"2f";"3e" **\dir{-};
"2f";"3f" **\dir{-};
"2e";"3f" **\dir{-};
"2e";"3g" **\dir{-};
\endxy
\]
\caption{The two (im)possible configurations of the edges in $K$}\label{f39}
\end{figure}

Now we return to the proof of Theorem  \ref{first}(\ref{c}), the case $p=2$, $q=3$.
Lemma \ref{quadrangle}, (a) and (b), in the preceding section, 
shows that there are at least four edges between $A$ and $B$. One case with
four edges is ruled out by having an isolated vertex in $B$, and another
by having two vertices in $B$ each with a single common neighbour in $A$,
these neighbours being the same. The remarks above rule out the complete
bipartite graph. The two cases that remain are shown in figure \ref{f39}.


Now we count some edges, along the same lines as earlier. The graph
induced on $\Omega f$ by $X$ omits three vertices,
and so, if $e$ denotes the number of edges of $X$, it has at most
$e-3k+3$ edges (with equality if the three vertices form a triangle).

We count the edges of $X f$ by seeing how many edges are identified
by $f$. All the edges within $K$ collapse to a single edge, so we lose
$3$ or $4$  edges. The other collapsed edges are
those with one end in $K$. The neighbours of a vertex in $A$ are mapped
injectively to $N(x)\setminus\{y\}$, and the neighbours of a vertex
of $B$ to $N(y)\setminus\{x\}$.

\paragraph{Case 1:} at most $k-2$ edges through $A$, and $(k-1)+(k-2)$
edges through $B$; thus at most $3k-2$ identified, at least $e-3k+2$ remain.

\paragraph{Case 2:} at most $k-3$ through $A$, and $(k-2)+(k-2)$ through
$B$; thus at most $3k-3$ identified, at least $e-3k+3$ remain.

\medskip

In Case 2, the fact that we have equality means that the neighbours outside
$K$ of the middle and bottom points of $B$ map to the same $k-2$ points of
$N(y)\setminus\{x\}$, so these sets are equal. Thus these two points of
$B$ have identical neighbour sets, contradicting primitivity (Lemma \ref{neigh}). So this case
cannot occur.

Now consider Case 1. Since the bounds differ by one, there are two subcases:

\subparagraph{Subcase 1A} The number of edges of $X f$ is $e-3k+2$,
which means that the neighbours outside $K$ of the two vertices in $A$ map
to the same $k-2$ neighbours of $x$, and so these two neighbour sets are equal.
In this case, the permutation interchanging the two vertices of $A$ and the
top and bottom vertices of $B$ is an automorphism, contradicting proper primitivity.

\subparagraph{Subcase 1B} The number of edges of $X f$ is $e-3k+3$, so
the bound is tight. In this case, the three vertices outside the image of $f$
must form a triangle. This cannot happen
if $f$ is an idempotent, since then the three vertices outside the image of $f$
lie in $K$, and cannot form a triangle.

The proof of Theorem \ref{first}(\ref{c}) is complete.\qed

\subsection{Proof of Theorem \ref{first}(\ref{d})}

The following result is the main observation underlying Theorem \ref{first}(\ref{d}).
\begin{lemma}
Let $f$ be a transformation on a finite set $\Omega$ and let $G$ be a group of permutations of $\Omega$. If there exists $g\in G$ such that $\rank(fgf)=\rank(f)$, then there exists an idempotent $e\in \langle f,G\rangle$ such that $e$ and $f$ have the same kernel.  
\end{lemma}
\begin{pf} Pick $g\in G$ such that $\rank(fgf)=\rank(f)$, that is, $|\Omega fgf|=|\Omega f|$.

Since $\langle f,G\rangle$ is finite it follows that there exists a smallest natural $n$ such that  $(fg)^{n}$ belongs to  $\{fg,(fg)^{2},(fg)^{3},\ldots, (fg)^{n-1}\}$.
Say $(fg)^{n}=(fg)^{m}$, with $m<n$.  We claim that  $(fg)^{n-m}$ is idempotent. 

First, observe that for all natural $l$ we have  $\rank((fg)^{l})=\rank(f)$. To see this 
suppose not and pick the smallest $i$ such that $\rank (fg)^{i}=\rank(f)$, but $\rank(fg)^{{i+1}}<\rank(f)$. 
Observe that  if $\rank(fg)^{i}=\rank (f)$, then we also have $\rank(f)=\rank(fg)^{i}=\rank((fg)^{i-1}f)$ (as a permutation in the end does not change the rank) and hence $\Omega (fg)^{i-1}f=\Omega f$. Now  $$|\Omega (fg)^{i+1}|=|\Omega (fg)^{i}fg|=|\Omega (fg)^{i}f|=|\Omega (fg)^{i-1}fgf|=|\Omega fgf|.$$ Thus $|\Omega f|=|\Omega (fg)^{i}|>|\Omega (fg)^{i+1}|=|\Omega fgf|,$ a contradiction with  the first sentence of this proof.  Thus $\rank(fg)^{l}=\rank (f)$, for all natural $l$.
In particular, $\rank (fg)^{m}=\rank (fg)^{n-m}=\rank(f)$. In addition, $\ker(fu)\supseteq \ker(f)$, for every transformation $u$. Thus $\ker(fg)^{m}$ and $\ker(fg)^{m-n}$ both contain $\ker(f)$. As their ranks are equal and we are dealing with finite sets, it follows that $\ker(fg)^{m}=\ker(fg)^{n-m}=\ker(f)$.

Now, $ (fg)^{m}=(fg)^{n}=(fg)^{m}(fg)^{n-m}$ implies that for all $x\in \Omega (fg)^{m}$ we have $x(fg)^{n-m}=x$. In addition, it is obvious that if $u=uv$, then $\Omega u \subseteq \Omega v$; as $(fg)^{m}=(fg)^{m}(fg)^{n-m}$ we get that $\Omega (fg)^{m}\subseteq \Omega (fg)^{n-m}$ and hence equality follows because $\rank(fg)^{m}=\rank(fg)^{n-m}$. It is proved that $\Omega(fg)^{m}=\Omega(fg)^{n-m}$ and we already know that $(fg)^{n-m}$ is the identity on $\Omega(fg)^{m}=\Omega(fg)^{n-m}$. It follows that $(fg)^{n-m}$ is idempotent. In addition, since we already proved that $\ker(fg)^{n-m} = \ker f$ the result follows.
\qed
\end{pf}

Now the proof of Theorem \ref{first}(\ref{d}) is immediate. If we have a group $G$ and a transformation $f$ under the hypothesis of the theorem, then, by the previous lemma, there exists an idempotent $e\in \langle f,G\rangle$ such that $\ker(f)=\ker (e)$. By  Theorem \ref{first}(\ref{c}) there exists a constant transformation $t$ such that $$t\in \langle e,G\rangle(\subseteq \langle f,G\rangle),$$ so that $G$ synchronizes $f$. 

\section{Proof of Theorem \ref{second}}\label{sma}

The aim of this section is to prove the following theorem.
\small*

Let $S$ be a transformation semigroup on a set of cardinality $n$ which
contains a transitive group $G$. If $f$ is an element of $S$ of minimum
rank $r$, then the image of $f$ is an $r$-clique in $\Gr(S)$, and the kernel
partition is a colouring with $r$ colours, and is uniform (so each part has
size $n/r$).

If $h$ is any element of $S$, then the rank of $h$ is at least $r$; we can
assume (replacing $f$ by $hf$ if necessary) that the kernel partition of $h$
refines that of $f$. We begin with a general result asserting that, if $G$
is primitive, then it is not possible for just one part of $f$ to be split
by the kernel partition of $h$.

\begin{theorem}
Let $S$ be a transformation semigroup containing a primitive group $G$, and
suppose that the minimum rank of an element of $S$ is $r$, where $r>1$. Then
there cannot be an element of $S$ with rank greater than $r$ whose kernel
partition has $r-1$ parts of size $n/r$.
\label{t5}
\end{theorem}
\begin{pf}
Let $f$ be an element of rank $r$ in $S$, and assume that $h$ is an element
whose kernel partition consists of $r-1$ parts of the kernel partition of $f$
and splits the remaining part into at least two.

Let $\Gamma=\Gr(S)$. We need one further observation about $\Gamma$. We know
that it has clique number $r$, so the independence number is at most $n/r$
(since $\Gamma$ is vertex-transitive). So, if $B$ is a part of the kernel
partition of $f$, and $v$ a vertex not in $B$, then $v$ has at least one
neighbour in $B$ (else $B\cup\{v\}$ would be an independent set).

Let $A_2,\ldots,A_r$ be the kernel classes of $h$ of size $n/r$, and
$A_{1,1},\ldots,A_{1,m}$ be the kernel classes into which the class $A_1$ of
$f$ is split. Let $a_j=A_jh$ and $a_{1,i}=A_{1,i}h$.

For $j,k>1$, there is an edge between $A_j$ and $A_k$; and
for $j>1$ and any $i$, there is an edge between $A_{1,i}$ and $A_j$. Since
$h$ is an endomorphism of $\Gamma$, there are edges between $a_j$ and $a_k$
for $j,k>1$, and between $a_{1,i}$ and $a_j$ for $j>1$ and all $i$. Thus
the subgraph on $\{a_{1,1},a_{1,2},a_2,\ldots,a_r\}$ is a complete
graph on $r+1$ vertices with an edge removed, contradicting Lemma~\ref{primgr}.\qed
\end{pf}
Let $S$ be a transformation semigroup containing a primitive group $G$, and
suppose that the minimum rank of an element of $S$ is $r$, where $r>1$. If $S$ contains a map $h$ of rank $r+1$, then either $\rank(hgh)=\rank(h)$, for all $g\in G$, or there exists $g\in G$ such that $\rank(hgh)=r$. The latter case cannot happen since  only two kernel blocks of $h$ collapse, and hence  $h$ would have $r-1$ kernel blocks of size $n/r$, in contradiction with the previous result. The former case implies that $h$ is uniform (by \cite{neumann}). But, for  $r>1$ it is impossible to have $(a_{1},\ldots,a_{r-1},b_{1},b_{2})$ and $(a_{1},\ldots,a_{r-1},b_{1}+b_{2})$  both uniform. We have proved the following corollary.

\begin{cor}\label{nonr+1}
Let $S$ be a transformation semigroup containing a primitive group $G$, and
suppose that the minimum rank of an element of $S$ is $r$, where $r>1$. Then
$S$ cannot contain an element of rank $r+1$.
\end{cor}
It is worth observing that this corollary immediately implies the result (proved by Rystsov) that the degree $n>1$ primitive groups on $\Omega$ synchronize  the rank $n-1$ transformations of $\Omega$.  In fact, if $n=2$, every rank $n-1$ map is already a constant and the result holds. If $n>2$, then a rank $n-1$ map $f$ cannot be uniform and hence, by \cite{neumann}, there exists $g\in G$ such that $\rank(fgf)<\rank(f)=n-1$. Thus let $n>n-1>r_k>\ldots >r_1$ be the possible ranks of the elements in $\langle G,f\rangle$. It is clear that for every $t\in \langle G,f\rangle$ and every $g\in G$ we have $\rank(tgf)\in\{\rank(t),\rank(t)-1\}$, since the kernel of $f$ does not allow to collapse more than two elements at once. This implies that $r_2=r_{1}+1$. Thus, by Corollary \ref{nonr+1}, $r_1$ cannot be larger that $1$ and hence $\langle G,f\rangle$ contains a constant, as claimed.

\paragraph{Proof of Theorem \ref{second}} As we have noted, the first part of the theorem
is due to P.M. Neumann \cite{neumann}. But that result can be  easily shown using the graph $\Gr(S)$. If the minimal rank
of an element in a transformation semigroup $S$ is $2$, then $\Gr(S)$ is
bipartite, and its automorphism group cannot be primitive (if $n>2$).

Since an element of minimal rank is uniform, we see that if $f$ is
non-uniform of rank~$3$ then the minimum rank in $\langle G,f\rangle$ is
either $1$ (so $G$ synchronizes $f$) or $2$ (so $G$ is imprimitive, since primitive groups synchronize rank $2$ maps).

Similarly, if $f$ is non-uniform of rank $4$ , then the minimum rank in $\langle G,f\rangle$
is $1$ (so $G$ synchronizes $f$), $2$ (so $G$ is imprimitive), or $3$ (so
the preceding corollary gives a contradiction).\qed

\section{Problems}\label{spro}
\setcounter{theorem}{0}

The major open question regarding the content of this paper is the following problem. 
\begin{prob}
Is it true that every  primitive group of permutations of a finite set $\Omega$  synchronizes a non-uniform transformation on $\Omega$?
\end{prob}

Assuming the previous question has an affirmative answer (as we believe), an 
 intermediate step in order to prove it would be to solve the following set of connected problems:
 
\begin{prob}
\begin{enumerate}
\item Remove the word
\emph{idempotent} in Theorem \ref{first}(\ref{c}).
 \item Extend Theorem \ref{first}(\ref{b}) to rank $n-3$.
 \item Prove that a primitive group synchronizes every non-uniform map of rank $5$. 
 \item\label{(d)} Prove that if $S=\langle f,G\rangle$ contains a map of minimal rank $r$, with $\rank(f)> r >1$,  there can be no map in $S$ of rank $r+2$. 
 
 Observe that if this conjecture is true, then the previous question is immediately also true. To see that, assume the conjecture holds and suppose that $f$ is a rank $5$ map.  By \cite{neumann}, if $r>1$, then $r\geq 3$; thus $r\in\{3,4\}$. Now the conjecture (if true) implies that $r+2\neq 5$, thus $r\neq 3$; on the other hand, Corollary \ref{nonr+1} implies that  $r\neq 4$. 
 This implies (modulo the conjecture) that a rank $5$ map $f$ either satisfies $\rank(fgf)=5$, for all $g\in G$, (and hence $f$ is uniform by \cite{neumann}), or $f$ is synchronized by any primitive group.  
 \end{enumerate}
\end{prob}

The next class of groups lies strictly between primitive and synchronizing. 

\begin{prob}
Is it possible to classify the primitive groups which synchronize
every rank $3$  map?
\end{prob}
		
		Note that there are primitive groups that do not synchronize a rank $3$ map (see the example immediately before Section \ref{trans}). And there are non-synchronizing groups which synchronize every rank $3$ map. Take for example $\pgl(2,7)$ of degree $28$; this group is non-synchronizing, but synchronizes every rank $3$  map	since $28$ is not divided by $3$.

Let $\Omega$ be a finite set of size $n$ and let $G$ be a non-synchronizing primitive group on $\Omega$. Let $s_1<s_2<\ldots <s_m$ be the  sequence of all the possible sizes of the $G$-regular partitions of $\Omega$ (that is, for every $s_i$ there exists a partition $P$ of $\Omega$, with $|P|=s_i$, and a set $S\subseteq \Omega$, such that $Sg$ is a section for $P$, for all $g\in G$). The depth of a group is 
$$
\mbox{d}(G):=\left\{
\begin{array}{ll}
s_2-s_1&\mbox{ if } m>1\\
\infty     &\mbox{ otherwise. }
\end{array}
\right.
$$
 Let $n$ be a natural number. If there exist non-synchronizing groups of degree $n$, then define the depth of $n$ as 
$$\mbox{d}(n):=\mbox{min}\{\mbox{d(G)}\mid G \mbox{ is a degree $n$ non-synchronizing primitive  group}\}.$$
\begin{prob}
\begin{enumerate}
\item Compute $\mbox{d}(n)$, for every natural number $n$ admitting degree $n$ primitive non-synchronizing groups. 
\item Let $G$ be a degree $n$ non-synchronizing primitive group with sequence $s_1<s_2<\ldots <s_m$ as above. Prove that $G$ synchronizes every  rank $s_1+e$ map (acting on the same set as $G$), where $$e\in \{1,\ldots,\mbox{d}(n)\}\cap\{1,\ldots ,n-1\}.$$ (Observe that Corollary \ref{nonr+1} already implies that this is true for $e=1$.)
\end{enumerate}

\end{prob}

There are very fast algorithms to decide if a given set of permutations generate a primitive group. 

\begin{prob}
Find an efficient algorithm to decide if a given set of permutations generate a synchronizing group.
\end{prob}

\begin{prob}\label{11}
Formulate and prove analogues of our results for semigroups of linear maps on a
vector space. Note that linear maps cannot be non-uniform, but we could ask for
linear analogues of results expressed in terms of rank, such as Theorems \ref{rys} and \ref{first}(\ref{b}).
\end{prob}

\begin{prob}
Solve the analogue of Problem \ref{11} for independence algebras (for definitions and fundamental results see \cite{Ar2,Ar4,Ar1,Ara1,Ara2,ArEdGi,arfo,Ar3,cameronSz,F1,F2,gould}). 
\end{prob}


\begin{thebibliography}{9}




\bibitem{Ar2}
J.~Ara\'ujo.
\newblock{Normal semigroups of endomorphisms of proper independence algebras are idempotent generated.}
\newblock{{\em  Proc. Edinburgh  Math. Soc.} \textbf{45} (2) (2002), 205--217. }

\bibitem{Ar4}
J.~Ara\'ujo.
\newblock{Generators for the semigroup of endomorphisms of an independence algebra.}
\newblock{{\em  Algebra Colloq.} \textbf{9} (4) (2002), 375--282. }

\bibitem{Ar1}
J.~Ara\'ujo.
\newblock{Idempotent generated endomorphisms of an independence algebra.}
\newblock{{\em Semigroup Forum} \textbf{67} (3) (2003), 464--467. }

\bibitem{Ara1}
J.~Ara\'ujo.
\newblock{Lifts for Semigroups of Monomorphisms of an Independence Algebra.}
\newblock{{\em  Colloquium Mathematicum.} \textbf{97} (2003), 277--284. }

\bibitem{Ara2}
J.~Ara\'ujo.
\newblock{Lifts for Semigroups of Endomorphisms of an Independence Algebra }
\newblock{{\em  Colloquium Mathematicum.} \textbf{106} (2006), 39--56. }


\bibitem{abc}
J. Ara\'ujo, W. Bentz and P.J. Cameron, Groups Synchronizing a Transformation of Non-Uniform Kernel.
\begin{verbatim}
http://arxiv.org/abs/1205.0682
\end{verbatim}

\bibitem{ArEdGi}
J.~Ara\'ujo, M. Edmundo and S. Givant.
\newblock{$v^*$-Algebras, Independence Algebras and Logic.}
\newblock{{\em International Journal of Algebra and Computation} \textbf{21} (7) (2011), 1237--1257. }


\bibitem{arfo}
J.\ Ara\'{u}jo and J.\ Fountain.
 \newblock{The Origins of Independence
Algebras}
\newblock \textit{Proceedings of the Workshop on Semigroups and Languages
(Lisbon 2002)}, World Scientific, (2004), 54--67

\bibitem{Ar3}
J.~Ara\'ujo and J.D. Mitchell.
\newblock{Relative ranks in the monoid of endomorphisms of an independence algebra.}
\newblock{{\em   Monatsh. Math.} \textbf{151} (1) (2007), 1--10. }


\bibitem{ArnoldSteinberg}
F. Arnold and B. Steinberg, Synchronizing groups and automata.
 \textit{  Theoret. Comput. Sci.} \textbf{359} (2006), no. 1-3, 101--110.

\bibitem{bamberg}
J. Bamberg, N. Gill, T.P. Hayes, H.A. Helfgott, A. Seress and P. Spiga,
Bounds on the diameter of cayley graphs of the symmetric
group.
\begin{verbatim}
http://arxiv.org/pdf/1205.1596.pdf
\end{verbatim}


\bibitem{CK}
P. J. Cameron and A. Kazanidis, 
Cores of symmetric graphs. \textit{J. Australian Math. Soc.} {\bf 85} (2008), 145--154.

\bibitem{cameronSz}
P.\ J.\ Cameron and C.\ Szab\'{o},
\newblock{Independence algebras},
\newblock{{\em J.\ London Math.\ Soc.},  \textbf{61} (2000),  321--334.}


\bibitem{DM}
J. D. Dixon and B. Mortimer, \textit{Permutation Groups}, Graduate Texts in Mathematics \textbf{163}, Springer, New York, 1996.
 
 
 \bibitem{F1}
J.~Fountain and A. Lewin.
\newblock{Products of idempotent endomorphisms of an independence algebra of finite rank.}
\newblock{{\em    Proc. Edinburgh Math. Soc.} \textbf{35} (2) (1992), 493--500. }

 \bibitem{F2}
J.~Fountain and A. Lewin.
\newblock{Products of idempotent endomorphisms of an independence algebra of infinite rank.}
\newblock{{\em   Math. Proc. Cambridge Philos. Soc.} \textbf{114} (2) (1993), 303--319. }



\bibitem{gould}
V.\ Gould,  
\newblock{Independence algebras.}
 \newblock{{\em Algebra
Universalis}  \textbf{33} (1995), 294--318.}

 
 \bibitem{helfgott}
H.A. Helfgott and A. Seress,
On the diameter of permutation groups, \emph{(to appear on the Annals of Math.)}
\begin{verbatim}
http://arxiv.org/abs/1109.3550
\end{verbatim}

 
\bibitem{levi}
I. Levi, R. McFadden, D. McAlister, 
Groups Associated with Finite Transformation Semigroups. \textit{Semigroup Forum} {\bf 61} (2000), 453--467.



\bibitem{neumann}
P. M. Neumann, Primitive permutation groups and their
section-regular partitions. \textit{Michigan Math. J.} {\bf 58} (2009), 309--322.

\bibitem{JEP}
J.-E. Pin, \v{C}ern\'{y}'s conjecture.
\begin{verbatim}
http://www.liafa.jussieu.fr/~jep/Problemes/Cerny.html
\end{verbatim}


\bibitem{rystsov}
I. Rystsov,
Quasioptimal bound for the length of reset words for regular automata.
\textit{Acta Cybernet.} {\bf 12} (1995), no. 2, 145--152. 


\bibitem{Tr}
A.N. Trahtman, Bibliography,  synchronization \@  TESTAS
\begin{verbatim}
http://www.cs.biu.ac.il/~trakht/syn.html
\end{verbatim}

\bibitem{Tr07}
A.N. Trahtman, \textit{The \v Cern\'y  Conjecture for Aperiodic Automata}, Discr. Math. \& Theoret. Comput. Sci. {\bf 9}(2007), (2)  3--10.

\bibitem{S-M}
Edward Szpilrajn-Marczewski, Sur deux propri\'et\'es des classes d'ensembles,
\textit{Fund. Math.} \textbf{33} (1945), 303-Ð307.

\end{thebibliography}
\end{document}